\documentclass[a4paper,twoside,12pt]{article}
\usepackage{graphicx,psfig,epsfig}
\usepackage{amsfonts}
\usepackage{amssymb}
\usepackage{amsmath}
\usepackage{subfigure}
\usepackage{setspace}
\usepackage{fancyheadings}
\usepackage[nosepfour,warning,np,debug,autolanguage]{numprint}
\usepackage[text={16cm,24cm},top=3cm,bottom=3cm,centering]{geometry}
\usepackage{algorithmicx}
\usepackage[ruled]{algorithm}
\usepackage{algc}

\newcommand\ASTART{\bigskip\noindent\begin{minipage}[b]{0.5\linewidth}}

\newcommand\AENDSKIP{\end{minipage}\bigskip}
\newcommand\AEND{\end{minipage}}
\usepackage{float}
\floatstyle{ruled}
\newfloat{Program}{thp}{loa}
\floatname{algorithm}{Algorithm}
\renewcommand{\maketitle}%
{
\noindent
\par\vspace*{3cm}\par
{\centering{\Large\textbf\MATCHtitle\par}
\par\vspace*{12pt}\par
{\textrm\MATCHauthor}\par\vspace*{12pt}\par}%
}
\renewenvironment{abstract}%
{\begin{center}\begin{minipage}[t]{1.0\textwidth}
\noindent{\centerline{\textbf{}}}
}%
{\end{minipage}\end{center}\par} 
\newenvironment{keywords}%
{\begin{center}\begin{minipage}[t]{1.0\textwidth}
\noindent\textbf{KEY WORDS$\colon$}}%
{\end{minipage}\end{center}\par}
 \makeatletter
 \renewcommand\section{\@startsection {section}{1}{\z@}%
                                   {-3.25ex \@plus -1ex \@minus -.2ex}%
                                   {2.3ex \@plus.2ex}%
                                   {\centering\large\bfseries\upshape}}
\makeatother 

\newcommand{\pfrac}[2]{\genfrac{(}{)}{}{}{#1}{#2}}
\usepackage{enumerate}
\usepackage[round,longnamesfirst]{natbib}
\usepackage[nosepfour,warning,np,debug,autolanguage]{numprint}
\newcommand{\MATCHtitle}%
{Multiblock Grid Generation for Simulations in Geological Formations\footnote[1]{Received\hspace*{2cm}; accepted\hspace*{2cm}. }}
\newcommand{\MATCHauthor}
{\textbf{Sanjay Kumar Khattri}\footnote[2]{Department of Mathematics, University of Bergen, Norway; e-mail$\colon$ sanjay@mi.uib.no; URL$\colon$ http://www.mi.uib.no/$\sim$sanjay.}}
\begin{document}
\begin{spacing}{1.0}
\maketitle
\begin{abstract}
\rule{460pt}{2pt}
{
Simulating fluid flow in geological formations requires mesh generation, lithology mapping to the cells, and computing geometric properties such as normal vectors and volume of cells. The purpose of this research work is to compute and process the geometrical information required for performing numerical simulations in geological formations. We present algebraic techniques, named Transfinite Interpolation, for mesh generation. Various transfinite interpolation techniques are derived from 1D projection operators. Many geological formations such as the Utsira formation \citep{Utsira_00,khattri_Utsira} and the  Sn{\o}hvit gas field \citep{snohvit} can be divided into layers or blocks based on the geometrical or lithological properties of the layers. We present the concept of block structured mesh generation for handling such formations.}
\rule{460pt}{2pt}
\end{abstract}
\section*{INTRODUCTION}
Simulation of fluid flow in geological formations, by numerical methods such as Finite Elements, Finite Volumes and Finite Differences, requires meshing of the geological formation into smaller elements called finite volumes or finite elements or cells depending on the numerical method \citep{khattri_Utsira,reservoir_comp_00,khattri_article_1,khattri_article_3,khattri_article_4}. These elements in three dimensions can be hexahedra, tetrahedra, prism and pyramid. In this paper, we focus only on hexahedral mesh generation. It is desirable that the part of the geological formation where solution shows nonlinear changes should be refined \citep{khattri_article_4,khattri_article_3,khattri_article_2}. Such a solution behaviour can occur due to lithological or geometrical properties of the formations \citep{khattri_article_4,khattri_article_3,khattri_article_2}. 

Many geological formations and reservoirs of interest can be divided into layers based on the geological characteristics such as faults and pinchouts or the lithological properties such as shale and sandstone. For example, the Utsira formation \citep{Utsira_00,khattri_Utsira} and the  Sn{\o}hvit gas field \citep{snohvit}. Each of these layers can be meshed into hexahedrals by the algebraic techniques independent of the other layers. In this way grid distribution and quality of mesh can be improved and controlled in each of the layers separately. This technique is called the multilayer or the multiblock approach. The concept of multiblock mesh generation is very useful for handling layered formations. Some of the advantages of this approach are 
\begin{enumerate}
\item Many geological formations can be realized by this concept.
\item It makes parallelization of a single phase problem straight forward. The multiblock/multilayer approach used as a domain decomposition concept allows the direct parallelization of both grid generation and flow codes on massively parallel systems. 
\item Grid density, distribution and quality can be controlled easily. It is desirable that in the areas of expected great nonlinear changes of solutions (around wells and material discontinuity) mesh should be refined. 
\item Controllability over the simulation. For example, the implementation of lithology and local optimization of mesh quality.
\item Though at the global level multilayer grids are unstructured in nature. Still at local level mesh can be expressed by logical numbering. Optimization of the quality of structured grids is easier. Instead of performing global mesh optimization, mesh can be optimized around critical locations such as wells. Structured grids can easily be made orthogonal at the boundaries and also almost orthogonal within the solution domain thus facilitating implementation of boundary conditions and also increase numerical accuracy. Discretization of partial differential equations on structured meshes is easier than on unstructured meshes. 
\item A structured grid produces a structured matrix and thus makes it easier to use sophisticated linear solvers. 
\end{enumerate}
Now let us discuss about algebraic method of grid generation.
\section*{Algebraic Method of Mesh Generation}
\begin{figure}
\hfill
\centering
\subfigure
{	
        \raisebox{0cm}{
    	\includegraphics[scale=1.0]{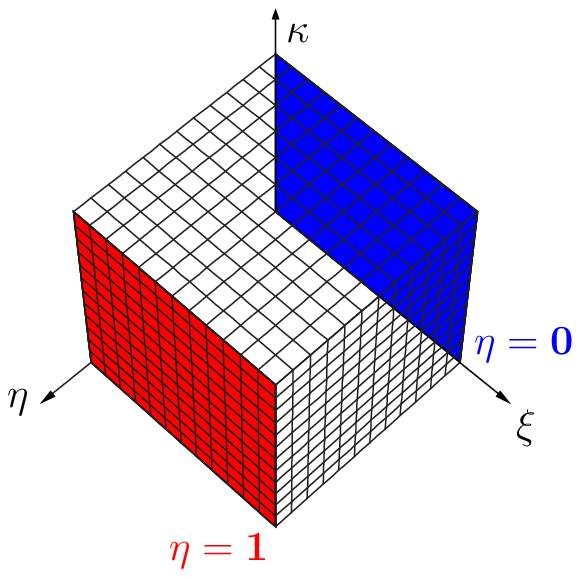}}
}
\centering
\subfigure
{	
	\raisebox{2.5cm}{
	\includegraphics[scale=1]{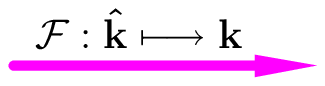}}
}
\hfill
\subfigure
{
    \includegraphics[scale=1.0]{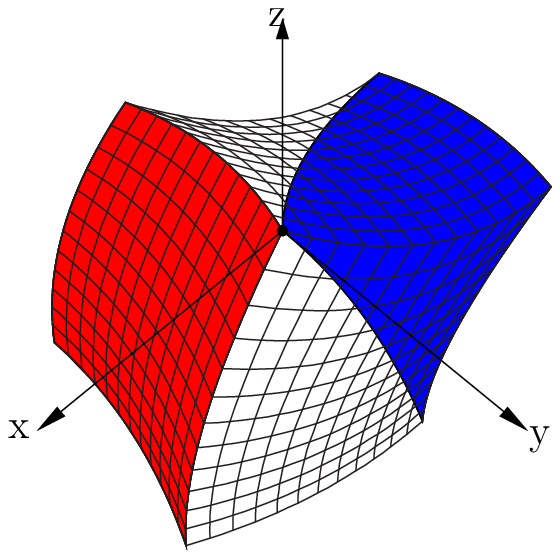}
}
\caption{Mapping a unit cube onto a physical domain.}
\label{fig:1}
\end{figure}

In the algebraic method of grid generation, we seek an algebraic mapping from a cube in computational or reference space to a physical space with the corresponding boundary surfaces \citep{MR1300634, MR1109550}. Transfinite interpolation (TFI) is such an algebraic mapping. TFI is also referred to as multivariate interpolation or Coons Patch. Figure \ref{fig:1} shows a mapping from a unit cube in the reference space onto a physical domain. Let the reference or computational space be defined by $\xi$, $\eta$ and $\kappa$ coordinates, and the physical space be defined by $x$, $y$ and $z$ coordinates. Suppose there exists a transformation or mapping, $\mathbf{r} = \mathbf{r}(\xi,\eta,\kappa)$, which maps the unit cube onto the interior of the physical domain, and this mapping maps the boundary surfaces of the cube to the corresponding boundary surfaces of the physical domain. Thus, $\eta=1$ surface of the cube is mapped to the $\mathbf{r}(1,\eta,\kappa)$ boundary surface of the physical domain. 
Transfinite interpolation is the boolean sum of univariate interpolations in each of the computational coordinates. Univariate interpolations are also referred to as one dimensional projection operators or projectors. Boolean sum of the projection operators are defined below. A univariate interpolation is an operator that vary only in one dimension or roughly speaking it is a function of only one reference coordinate. A univariate interpolation can be linear, quadratic and cubic. Any univariate interpolation can be applied in a coordinate direction. Generally a higher order interpolation operator is desired in flow direction. TFI is composed of 1D projection operators, let us first define some one dimensional projection operators.
\section*{One Dimensional Projection Operators}
A 1D projection operator or projector can be defined in many ways depending upon the available information. For example, a linear projector can be formed from two surfaces; a Hermite projector can be formed from two surfaces and directional derivatives at these surfaces; a Lagrangian projector can be defined from two boundary surfaces and internal surfaces.  

Let the reference space be defined by $\xi$, $\eta$ and $\kappa$ coordinates ($\xi\in[0,1]$, $\eta\in[0,1]$ and $\kappa\in[0,1]$). Suppose there exists a transformation $\mathbf{r}(\xi,\eta,\kappa)$ from a unit cube in the reference space onto a physical domain. That is $\mathbf{r}\colon \hat{k}\longmapsto{k}$. Let the physical space be defined by six boundary surfaces. A $\xi$ surface in the physical space is a surface on which value of $\xi$ is constant. Thus, two $\xi$ boundary surfaces are $\mathbf{r}(0,\eta,\kappa)$ and $\mathbf{r}(1,\eta,\kappa)$. Similarly, two $\eta$ and $\kappa$ boundary surfaces are given as $\mathbf{r}(\xi,0,\kappa)$, $\mathbf{r}(\xi,1,\kappa)$ and $\mathbf{r}(\xi,\eta,0)$, $\mathbf{r}(\xi,\eta,1)$, respectively. From these six boundary surfaces, the following 1D projection operators are defined
\begin{alignat}{3}
{\mathbf{P}_{\xi}} &\overset{\textbf{def}}{=}(1- \xi) \, {\mathbf {r}(0,\eta,\kappa)} + \xi \, {\mathbf {r}(1,\eta,\kappa)}\enspace, \label{eq1}  \\
{\mathbf {P}_{\eta}}  &\overset{\textbf{def}}{=}(1- \eta) \, {\mathbf {r}(\xi,0,\kappa)}  + \eta \,  {\mathbf {r}(\xi,1,\kappa)}\enspace,\label{eq2}\\
{\mathbf {P}_{\kappa}} &\overset{\textbf{def}}{=}(1- \kappa) \, {\mathbf {r}(\xi,\eta,0)} + \kappa \, {\mathbf {r}(\xi,\eta,1)}\enspace.\label{eq3}
\end{alignat}
The projectors $\mathbf{P_{\xi}}$, $\mathbf{P_{\eta}}$ and $\mathbf{P_{\kappa}}$ are 1D projection operators and they are functions of the coordinates ($\xi,\eta,\kappa$). 
The projection operators defined by equations \eqref{eq1}, \eqref{eq2} and \eqref{eq3} are linear in $\xi$, $\eta$ and $\kappa$ coordinates. It can be notice that the operators are defined from two surfaces of a particular kind. For example, $\mathbf{P}_{\xi}$ is defined from two $\xi$ boundary surfaces in the physical space $\mathbf{r}(0,\eta,\kappa)$ and $\mathbf{r}(1,\eta,\kappa)$. 

If in addition to the boundary surfaces we also know the internal surfaces of a domain then a projection operator can also be defined from more than two surfaces of a kind. For example, if there are $n+1$ surfaces of $\xi$ type ($n-1$ internal curves and $2$ boundary surfaces) then $\mathbf{P}_{\xi}$ projection operator can be defined as
\begin{align}
\mathbf{P}_{\xi} &\overset{\textbf{def}}{=} \sum^{n}_{j=0} \beta_{j}(\xi) \, \mathbf{r}(\xi_j,\eta,\kappa)\label{eq:lagprojop}\enspace, \\
\intertext{\citep{lagrange1,lagrange2}. 
Here, $j=0$ and $j=n$ are the boundary surfaces while $j=1,\ldots,n-1$ are the internal surfaces, and $\beta_{j}$ is the Lagrangian weighting factor. The Lagrangian weighting factor is given as follows}
\beta_{j}(\xi) &= \prod_{{i=0,\, i \ne j}}^n \pfrac{\xi-\xi_i}{\xi_j-\xi_i}\label{eq:weight_largange}\enspace.
\\
\intertext{It can be notice that the weighting factor $\beta_j(\xi)$ is an order $n$ polynomial having zeros at all of the surfaces except the $jth$ surface. The Lagrangian weighting factor satisfies the following}
\beta_{j}(\xi_{i}) = 
\begin{cases}
1 \quad & \text{if $i = j$}\enspace,\\
0 \quad & \text{if $i \ne j$}\enspace,
\end{cases}\quad & \text{and} \quad \sum^{n}_{j=0}\beta_j = 1.0\enspace.
\end{align}
Now let us express the Lagrangian projection operator in another form. The numerator in the Lagrange weighting factor \eqref{eq:weight_largange} can be written as 
\begin{equation}
\dfrac{\left[(\xi-\xi_0)\,(\xi-\xi_1)\,\cdots\,(\xi-\xi_n)\right]}{(\xi-\xi_j)} = \dfrac{\Omega}{(\xi-\xi_j)}\label{eq:num1}\enspace,
\end{equation}
\citep[see][]{lagrange1}. Let us define the barycentric weights as \citep{lagrange1}
\begin{equation}
\omega_j = \dfrac{1}{\prod_{i=0,\,i\ne j}^n(\xi_j-\xi_i)}\label{eq:weight2}\enspace.
\end{equation}
Using equations \eqref{eq:num1} and \eqref{eq:weight2}, the Lagrangian projection operator \eqref{eq:lagprojop} can also be written as \citep{lagrange1}
\begin{equation}
\mathbf{P}_{\xi} \overset{\textbf{def}}{=} \Omega \sum^{n}_{j=0} \dfrac{\omega_j}{\xi-\xi_j} \, \mathbf{r}(\xi_j,\eta)\label{eq:lagprojop1}\enspace.
\end{equation}
In the grid generation literature, the equation \eqref{eq:lagprojop} is used but the new form \eqref{eq:lagprojop1} is computationally more efficient \citep[cf.][]{lagrange1}.
Similarly, if in addition to the boundary surfaces we are also given the derivatives (direction vectors) on these boundary surfaces then we can define the Hermite interpolation operators. For example, if we are given two $\xi$ surfaces$\colon$ $\mathbf{r}(0,\eta,\kappa)$ and $\mathbf{r}(1,\eta,\kappa)$, and let the direction vectors on these surfaces be $\mathbf{r{'}}(0,\eta,\kappa)$ and $\mathbf{r}'(0,\eta,\kappa)$, respectively. Then, the 1D Hermite projection operator can be defined as
\begin{alignat}{2}
\mathbf{P}_{\xi} \overset{\textbf{def}}{=} (2\,{\xi}^{3}-3\,{\xi}^{2}+1) \, \mathbf{r}(0,\eta,\kappa) &{}+ (-2\,{\xi}^{3}+3\,{\xi}^{2})\, \mathbf{r}(1,\eta,\kappa) \nonumber\\ 
&{}+ ({\xi}^{3}-2\,{\xi}^{2}+\xi) \, \mathbf{r{'}}(0,\eta,\kappa)+({\xi}^{3}-{\xi}^{2})\, \mathbf{r{'}}(1,\eta,\kappa)\enspace.\label{eq4} 
\end{alignat}
Hermite projectors are easy to implement and are powerful tools for grid generation. Grid lines can be made orthogonal by the proper choice of direction vectors. This may help in accurate modelling of boundary conditions. 
\begin{figure}
\begin{minipage}[b]{0.45\linewidth} 
\centering
\includegraphics[scale=1.0]{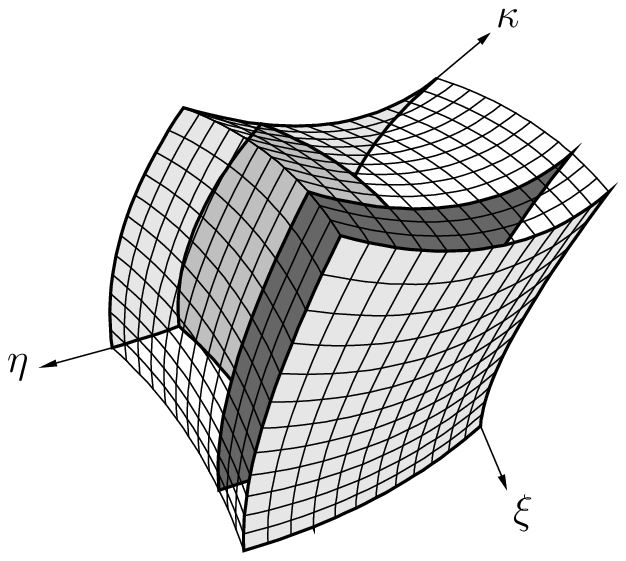} 
\caption{A 3D physical domain containing 3 $\xi$, 3 $\eta$ and 2 $\kappa$ surfaces.}
\label{fig:physical}
\end{minipage}
\hspace{0.5cm} 
\begin{minipage}[b]{0.45\linewidth}
\centering
\includegraphics[scale=1.0]{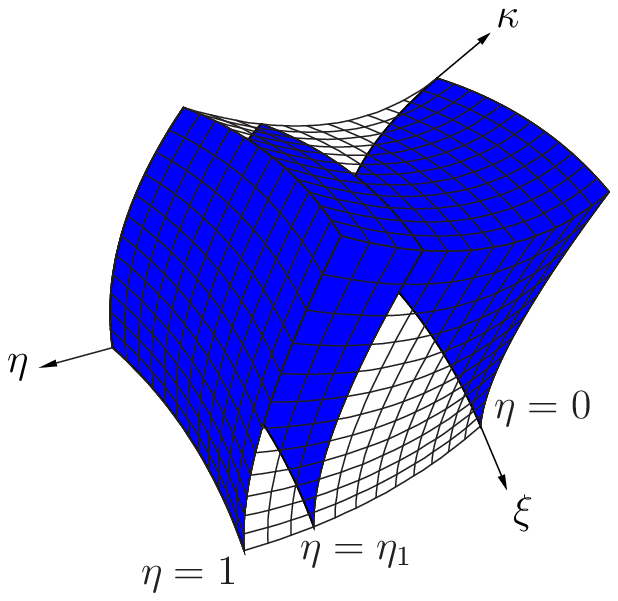}
\caption{A 3D physical domain containing 2 $\xi$, 3 $\eta$ and 2 $\kappa$ surfaces.}
\label{fig:physical_1}
\end{minipage}
\end{figure}

Figure \ref{fig:physical} shows a physical domain containing three $\xi$, three $\eta$ and two $\kappa$ surfaces. Since the domain contains three $\xi$ surfaces, three $\eta$ surfaces and two $\kappa$ surfaces thus we can define a Lagrangian $\mathbf{P}_\xi$ operator, a Lagrangian $\mathbf{P}_\eta$ operator and a linear $\mathbf{P}_\kappa$ operator. Figure \ref{fig:physical_1} shows another physical domain with two $\xi$, three $\eta$ ($\mathbf{r}(\xi,0,\kappa)$, $\mathbf{r}(\xi,\eta_1,\kappa)$ and $\mathbf{r}(\xi,1,\kappa)$) and two $\kappa$ surfaces. For this domain, a linear $\mathbf{P}_\xi$, a Lagrangian $\mathbf{P}_\eta$ and a linear $\mathbf{P}_\kappa$ operators can be defined. For this domain, the Lagrangian $\mathbf{P}_\eta$ operator is given as
\begin{equation}
\mathbf{P}_\eta = \Omega 
\left[\left(\dfrac{\omega_0}{\eta-0}\right)\mathbf{r}(\xi,0,\kappa) + \left(\dfrac{\omega_1}{\eta-\eta_1}\right)\mathbf{r}(\xi,\eta_1,\kappa) + 
\left(\dfrac{\omega_2}{\eta-1}\right)\mathbf{r}(\xi,1,\kappa) \right]\enspace,
\end{equation}
where $\Omega$ is given as,
$$\Omega = \eta\,(\eta-\eta_1)\,(\xi-\xi_3)\enspace,$$
{and $\omega_0$, $\omega_1$, $\omega_2$ and $\omega_3$ are given as,}
\begin{alignat}{4}
\omega_0 &= \dfrac{1}{\eta_1}\enspace, & \qquad 
\omega_1 &= \dfrac{1}{(-\eta_1)\,(1-\eta_1)}\enspace,\qquad
\omega_2 &= \dfrac{1}{(-1)\,(\eta_1-1)}\enspace. & \qquad 
\end{alignat}

Now let us study two important and useful properties of projection operators. These properties are called tensor product and boolean sum of projection operators.
\section*{Properties of Projection Operators}
This section presents two important properties of projection operators.
\subsection{Tensor Product} Tensor product $\mathbf{P}_{\xi\circ\eta}$ of the projection operators $\mathbf{P}_{\xi}$ and $\mathbf{P}_{\eta}$ is defined as follows
\begin{alignat}{4}
\mathbf{P}_{\xi\circ\eta} &\overset{\textbf{def}}{=} \mathbf{P}_{\xi} \circ  \mathbf{P}_{\eta} = (1-\xi)\,[\mathbf{P}_{\eta}]_{\xi = 0} + \xi\,[\mathbf{P}_{\eta}]_{\xi = 1}\enspace. \label{prop1} \\
\intertext{Here, $\mathbf{P}_{\xi}$ is assumed to be linear projection operator as defined by the equation \eqref{eq1}. It is clear from equation \eqref{prop1} that $\mathbf{P}_\xi$ is a projection operator. That is $\mathbf{P}_{\xi\circ\xi}$ = $\mathbf{P}_{\xi}$. If $\mathbf{P}_\xi$ is Lagrangian projection operator then the tensor product is defined as}
\mathbf{P}_{\xi\circ\eta} &\overset{\textbf{def}}{=} \mathbf{P}_{\xi} \circ  \mathbf{P}_{\eta} = 
\sum_{j=0}^n\beta_j(\xi)\,\left[\mathbf{P}_\eta\right]_{\xi=\xi_j}\enspace.	
\end{alignat}
Tensor product of two projection operators is also a projection operator ($\mathbf{P_{\xi\circ\eta}}$ is a projection operator). Since tensor product is also a projection operator, it is commutative in nature. That is $\mathbf{P}_{\xi\circ\eta}$ = $\mathbf{P}_{\eta\circ\xi}$. Similarly tensor products can be defined for an arbitrary number of projection operators. For example, the tensor product of three projection operators is defined as follows
\begin{equation}
\mathbf{P}_{\xi\circ\eta\circ\kappa} \overset{\textbf{def}}{=} \mathbf{P}_{\xi} \circ  (\mathbf{P}_{\eta} \circ \mathbf{P}_{\kappa})=(1-\xi)\,[\mathbf{P}_{\eta\circ\kappa}]_{\xi = 0} + \xi\,[\mathbf{P}_{\eta\circ\kappa}]_{\xi = 1}\enspace.
\end{equation}
In the above equation, the projection operator $\mathbf{P}_\xi$ is linear.
\subsection{Boolean Sum} Boolean sum of two projection operators is a also a projection operator and it is defined as follows
\begin{equation}
\mathbf{P}_{\xi\oplus\eta}\overset{\textbf{def}}{=} \mathbf{P}_{\xi} \oplus \mathbf{P_{\eta}} = \mathbf{P}_{\xi} +  \mathbf{P}_{\eta} - \mathbf{P}_{\xi\circ\eta}\enspace.
\end{equation}
Here, $\mathbf{P}_{\xi\circ\eta}$ is the tensor product of the $\mathbf{P}_{\xi}$ and  $\mathbf{P}_{\eta}$ projection operators. Boolean sum is commutative in nature. That is $\mathbf{P}_{\xi}$ $\oplus$ $\mathbf{P}_{\eta}$ = $\mathbf{P}_{\eta}$ $\oplus$ $\mathbf{P}_{\xi}$. Since boolean sum is also a projection operator thus it follows the projection property. That is $\mathbf{P}_{\xi\oplus\xi}$ = $\mathbf{P}_{\xi}$. Similarly, the boolean sum can also be defined for an arbitrary number of projection operators. The boolean sum of three projectors is defined by using the fact that boolean sum and tensor product of two projection operators are also projection operators. Thus, the boolean sum of $\mathbf{P}_{\xi}$, $\mathbf{P}_{\eta}$ and $\mathbf{P}_{\kappa}$ operators is given as
\begin{alignat}{2}
\mathbf{P}_{\xi\oplus\eta\oplus\kappa} &=\mathbf{P_{\xi}} \oplus \mathbf{P_{\eta}} \oplus \mathbf{P_{\kappa}}\enspace,\nonumber\\
& = \mathbf{P}_{\xi} \oplus\left( \mathbf{P}_{\eta} \oplus \mathbf{P}_{\kappa}\right)\nonumber\enspace,\\
&= \mathbf{P}_{\xi} \oplus\left( \mathbf{P}_{\eta} + \mathbf{P}_{\kappa} - \mathbf{P}_{\eta\circ\kappa}\right)\enspace,\nonumber\\
&= \mathbf{P}_{\xi} \oplus \mathbf{P}_{\eta} + \mathbf{P}_{\xi} \oplus \mathbf{P}_{\kappa} - \mathbf{P}_{\xi} \oplus \mathbf{P}_{\eta\circ\kappa} \enspace,\nonumber\\
&= \mathbf{P}_{\xi} + \mathbf{P}_{\eta} - \mathbf{P}_{\xi\circ\eta} + \mathbf{P}_{\xi} + \mathbf{P}_{\kappa} - \mathbf{P}_{\xi\circ\kappa} - \mathbf{P}_{\xi} -\mathbf{P}_{\eta\circ\kappa} + \mathbf{P}_{\xi\circ\eta\circ\kappa} \enspace.\nonumber\\
\intertext{Thus,}
\mathbf{P}_{\xi\oplus\eta\oplus\kappa}&= \mathbf{P}_{\xi} + \mathbf{P}_{\eta} +\mathbf{P}_{\kappa} - \mathbf{P}_{\xi\circ\eta} - \mathbf{P}_{\xi\circ\kappa} -\mathbf{P}_{\eta\circ\kappa} + \mathbf{P}_{\xi\circ\eta\circ\kappa} \enspace. \label{eq:bool_sum_00}
\end{alignat}
Here, $\mathbf{P}_{\xi\circ\eta}$ denotes the tensor product of $\mathbf{P}_{\xi}$ and $\mathbf{P}_{\eta}$ projection operators and $\mathbf{P}_{\xi\circ\eta\circ\kappa}$ denotes the tensor product of $\mathbf{P}_{\xi}$, $\mathbf{P}_{\eta}$ and $\mathbf{P}_{\kappa}$ projection operators.
\section*{Transfinite Interpolation}
Boolean sum of projection operators is the basis for Transfinite Interpolation. TFI are extensible used for algebraic grid generation. Since, 1D projection operators comes in many flavours such as the Lagrangian and the Hermite thus TFI can be defined by many different expressions depending upon which 1D projection operators are used. Linear Transfinite Interpolation creates a grid in 3D using surfaces that define the boundaries. Quality of the generated grid strongly depends on the parametrizations of the boundary curves. In its simplest form this mapping blends two given surfaces to create a grid in the region bounded by the surfaces or curves. Linear Transfinite Interpolation mapping is defined from six surfaces. Transfinite Interpolation mapping will only give a reasonable grid if the surfaces that define the boundary match at the edges, and the surfaces are parametrized in the same direction otherwise grid lines could cross each other. We are using the equation \eqref{eq:bool_sum_00} for mesh generation. Thus, the position vector in the physical space is given as
\begin{equation}
\mathbf r(\xi,\eta,\kappa) = \mathbf{P}_{\xi\oplus\eta\oplus\kappa} = \mathbf{P}_{\xi} \oplus \mathbf{P}_{\eta} \oplus \mathbf{P}_{\kappa} \enspace.
\label{eq:tranfinite_intp}
\end{equation}
Let the geological formation be defined by the six boundary surfaces $\mathbf{r}(0,\eta,\kappa)$, $\mathbf{r}(1,\eta,\kappa)$, $\mathbf{r}(\xi,0,\kappa)$, $\mathbf{r}(\xi,1,\kappa)$, $\mathbf{r}(\xi,\eta,0)$ and $\mathbf{r}(\xi,\eta,1)$. Thus, from these six boundary surfaces the linear projection operators can be defined. Let us divide the reference unit cube into $nx$ subdivisions in the $\xi$ coordinate direction, $ny$ subdivisions in the $\eta$ coordinate directions, and $nz$ subdivisions in the $\kappa$ coordinate direction. Thus for this mesh 
\begin{alignat}{3}
\text{Number of nodes} &=& &nx\times{ny}\times{nz},\nonumber\\
\text{Number of cells} &=&  &(nx+1)\times{(ny+1)}\times{(nz+1)},\nonumber\\
\text{Number of surfaces} &=& &nx\times ny \times (nz+1) + nx \times nz \times (ny+1) +  ny \times nz \times (nx+1). \nonumber   
\end{alignat}
A simple routine for generating mesh in the geological formation is given as
\begin{algorithm}
\caption{Grid generation in a block or layer.}
\begin{algorithmic}[1]
\For{\mathrm{ix=0} }{\mathrm {ix < nx+1} }{\mathrm {ix^{++}} }{\bf{\{}} {\Comment\textbf{Moving in the $\xi$ direction}}
\For{\mathrm{iy=0} }{\mathrm{iy < ny+1} }{\mathrm{iy^{++}} }{\bf{\{}} {\Comment\textbf{Moving in the $\eta$ direction}}
\For{\mathrm{iz=0}}{\mathrm{iz < nz+1} }{\mathrm{iz^{++}} }{\bf{\{}} {\Comment\textbf{Moving in the $\kappa$ direction}} \\
\quad \quad \quad \quad $\mathrm{i := ix + (nx+1) \times iy + (ny+1) \times iz}$; {\Comment\textbf{Node number}} \\
\quad \quad \quad \quad $\mathrm{\xi_1 := {ix}/{nx};\quad \eta_1 := {iy}/{nx};\quad \kappa_1 := {iz}/{nx}};$ 
{\Comment\textbf{Gridding of Unit Cube}}\\
\quad \quad \quad \quad $\mathrm{\mathbf{r}(ix,iy,iz):= {[\mathbf{P}_{\xi\oplus\eta\oplus\kappa}]}_{\xi=\xi_1,\eta=\eta_1,\kappa=\kappa_1}}$ {\Comment\textbf{Position in the Physical Space}}\\
\quad \quad \quad {\bf{\}}}\\
\quad \,\, {\bf{\}}}\\
{\bf{\}}}
\end{algorithmic}
\end{algorithm}
\section*{Computing Geometric Properties}
\label{sec:comp_geom}
\begin{figure}
\begin{minipage}[b]{0.31\linewidth} 
\centering
\includegraphics[scale=0.40]{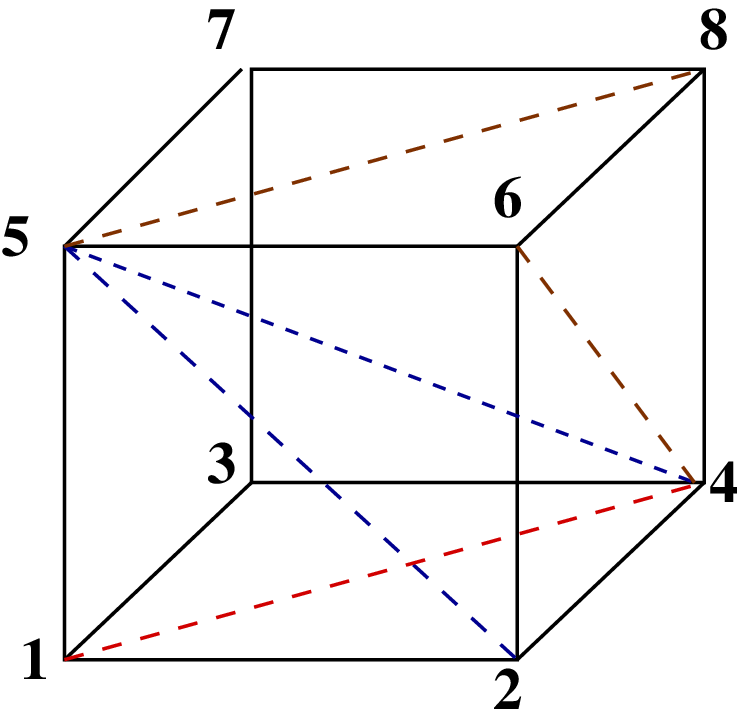} 
\caption{Division of a hexahedra.}
\label{fig:divide_hex}
\end{minipage}
\begin{minipage}[b]{0.31\linewidth}
\centering
\includegraphics[scale=0.50]{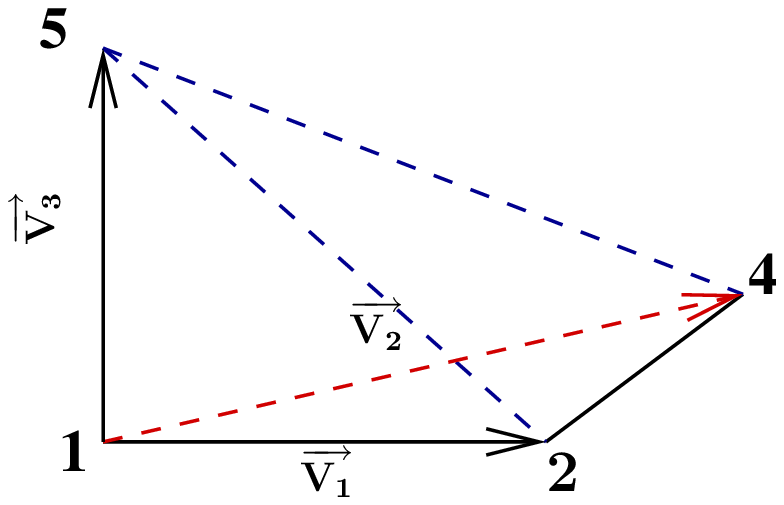} 
\caption{Volume of the tetrahedra 1245.}
\label{fig:lonely_tetrahedra} 
\end{minipage}
\begin{minipage}[b]{0.31\linewidth}
\centering
\includegraphics[scale=0.5]{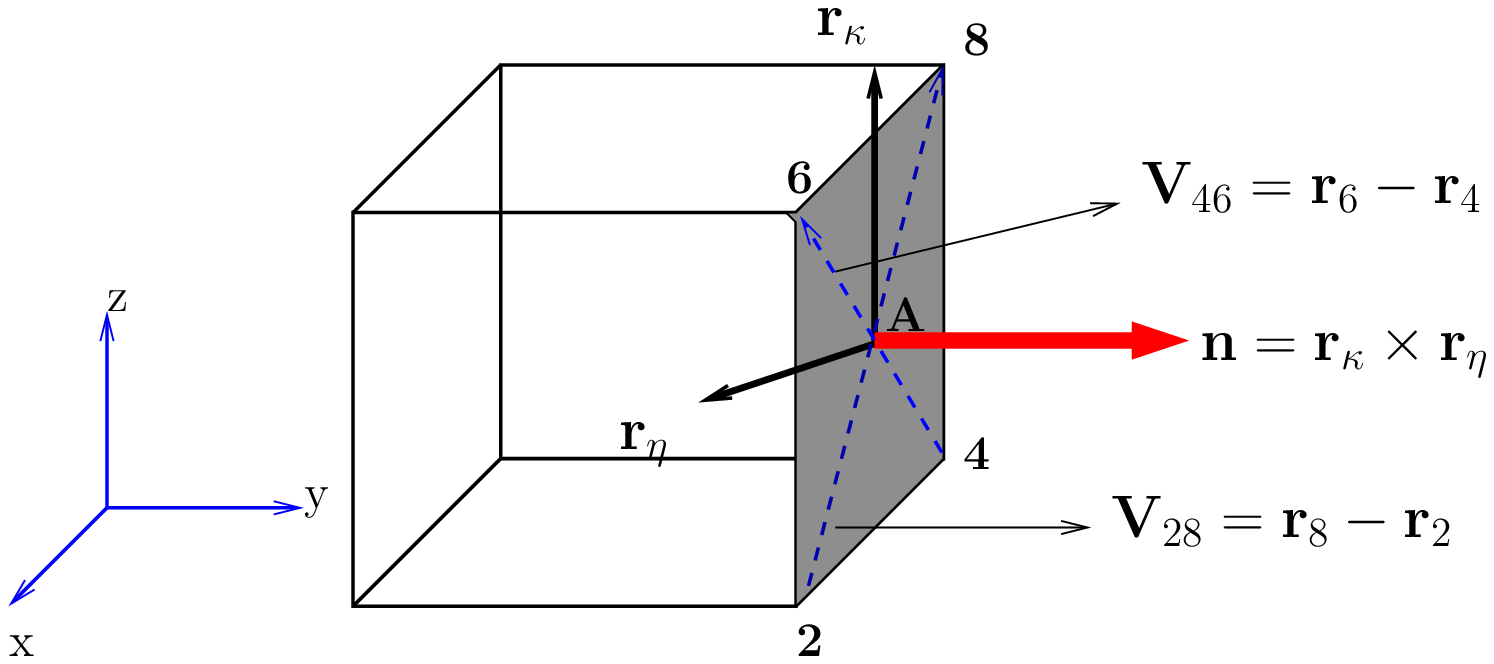}
\caption{Normal vector on the surface 1234.}
\label{fig:normal_vec}
\end{minipage}
\end{figure}
Let us consider the steady state pressure equation of a single phase flowing in a porous medium 
\begin{alignat}{2}
\label{elliptic1}
-\,\text{div}\,(\boldsymbol{K}\,\text{grad}\,p) &= f \enspace{.}
\end{alignat}
In porous media flow, the unknown function $p=p(x,y)$ represents the pressure of a single phase, $\boldsymbol{K}$ is the permeability or  hydraulic conductivity of the porous medium, and the velocity $\mathbf{u}$ of the phase is given by the Darcy law as$\colon$ $\mathbf{u}=-\boldsymbol{K}\,{\text{grad}\,{p}}$. For solving partial differential equations (PDEs) in geological formations by numerical methods such as Finite Volumes, the domain is divided into smaller elements. The process of dividing geological formations into smaller elements is referred to as meshing of the domains or geological formations, and the elements are called finite volumes or cells. 
Integrating equation \eqref{elliptic1} over one of the finite volumes with volume {\bf{Vol}} and boundary {$\partial\textbf{Vol}$}, and using the Gauss divergence theorem leads to 
\begin{alignat}{3}
-\int_{\partial\textbf{Vol}}{\boldsymbol{K}\,{\nabla{p}}}\cdot{\mathbf{\hat{{n}}}} &= \int_{\textbf{Vol}}{f}\enspace,  \label{gauss-D}\\
\intertext{where ${\mathbf{\hat{n}}}$ is the outward unit normal on the boundary ${\partial\textbf{Vol}}$ of the finite volume \textbf{Vol}. Let us assume that finite volumes are hexahedras. Boundary of these finite volumes consists of six surfaces $\partial\textbf{Vol}_i$. The above equation can be written as}
-\sum_{i=1}^6{\int_{\partial\textbf{Vol}_{i}}{\boldsymbol{K}\,{\nabla{p}}}\cdot{\mathbf{\hat{n}}}} &= \int_{\textbf{Vol}}{f}\enspace, \label{discrete_3}
\end{alignat}
the term $-{\int_{{\partial\textbf{Vol}_{i}}}{\boldsymbol{K}\,{\nabla{p}}}\cdot{\mathbf{\hat{n}}}}$ is referred to as the flux or the Darcy flux through the surface ${\partial\textbf{Vol}_{i}}$. The term $\int_{\textbf{Vol}}{f}$ can be approximated as value of the function $f$ at the center of the hexahedra times the volume of the hexahedra. Thus, converting a partial differential equation into an algebraic equation requires volume of hexahedra and normal vectors on the surfaces of the hexahedra. Now, let us present a method for computing volume of the hexahedra.

Figure \ref{fig:divide_hex} shows a hexahedra {\bf{12345678}}. Let the position vector of the vertix ${{i}}$ be $\mathbf{r}_i$ with $i={\bf{1,\ldots,8}}$. This hexahedra can be divided into two prisms {\bf{124568}} and {\bf{134578}}. Each of these prisms can divided into three tetrahedras.  The Figure \ref{fig:divide_hex} shows the division of the prisms {\bf{124568}} into three tetrahedras {\bf{1245}}, {\bf{2456}} and {\bf{4568}}. Thus, a hexahedra can divided into six tetrahedras, and the volume of the hexahedra can be computed by summing the volume of the six tetrahedras. Figure \ref{fig:lonely_tetrahedra} presents the tetrahedra {\bf{1245}}. The vectors $\overrightarrow{\mathbf{V}}_1$, $\overrightarrow{\mathbf{V}}_2$ and $\overrightarrow{\mathbf{V}}_3$ are meeting at the vertix 1 of the tetrahedra. The vectors $\overrightarrow{\mathbf{V}}_1$, $\overrightarrow{\mathbf{V}}_2$ and $\overrightarrow{\mathbf{V}}_3$ are given as $\overrightarrow{\mathbf{V}}_1$ = $\mathbf{r}_2-\mathbf{r}_1$, $\overrightarrow{\mathbf{V}}_2$ = $\mathbf{r}_4-\mathbf{r}_1$ and $\overrightarrow{\mathbf{V}}_3$ = $\mathbf{r}_5-\mathbf{r}_1$, respectively. The volume of the tetrahedra {\bf{1245}} is given as 
\begin{equation}
\textbf{Vol}_{\bf{1245}} = \dfrac{1}{6}\,\vert{\overrightarrow{\mathbf{V}}_1}\cdot(\overrightarrow{\mathbf{V}}_2\times\overrightarrow{\mathbf{V}}_3)\vert\enspace.
\end{equation}
Now, we are going to see two techniques for computing normal vectors on the surface of hexahedra. 

For the surface $\bf{2468}$, see Figure \ref{fig:normal_vec}. The diagonal vectors $\mathbf{V}_{28}$ and $\mathbf{V}_{46}$ of the quadrilateral surface {\bf{2486}} of the hexahedra are given as $\mathbf{V}_{28} = \mathbf{r}_8-\mathbf{r}_2$ and $\mathbf{V}_{46} = \mathbf{r}_{6}-\mathbf{r}_4$. The normal vector on the quadrilateral surface is given as the cross product of these two diagonal vectors. That is $\mathbf{n}=\mathbf{V}_{28}\times\mathbf{V}_{46}$.

The position vector of a point in the physical space (geological formation) is given by the expression \eqref{eq:tranfinite_intp}, and this expression is a function of the coordinates $\xi$, $\eta$ and $\kappa$. Differentiating this expression with respect to a particular coordinate will give us a vector pointing in that coordinate direction. This vector is called the covariant vector. Figure \ref{fig:normal_vec} presents two covariant vectors $\mathbf{r}_\eta$ and $\mathbf{r}_\kappa$. Differentiating the expression \eqref{eq:tranfinite_intp} with respect to $\eta$ results
\begin{equation}
\mathbf{r}_\eta = \dfrac{\partial{\mathbf{P}_\eta}}{\partial\eta} - \dfrac{\partial{\mathbf{P}_{\xi\circ\eta}}}{\partial\eta}  - \dfrac{\partial{\mathbf{P}_{\eta\circ\kappa}}}{\partial\eta} +
\dfrac{\partial{\mathbf{P}_{\xi\circ\eta\circ\kappa}}}{\partial\eta}\enspace.
\end{equation}
Since, $\mathbf{P}_\xi$ and $\mathbf{P}_\kappa$ are not functions of $\eta$ so their differentiation with respect to $\eta$ will vanish. Similarly, the covariant vector $\mathbf{r}_\kappa$ can be determined. Cross product of these two covariant vectors will provide the normal vector on the surface. 
\section*{Example}
\begin{figure}
\centering
\includegraphics[scale=0.5]{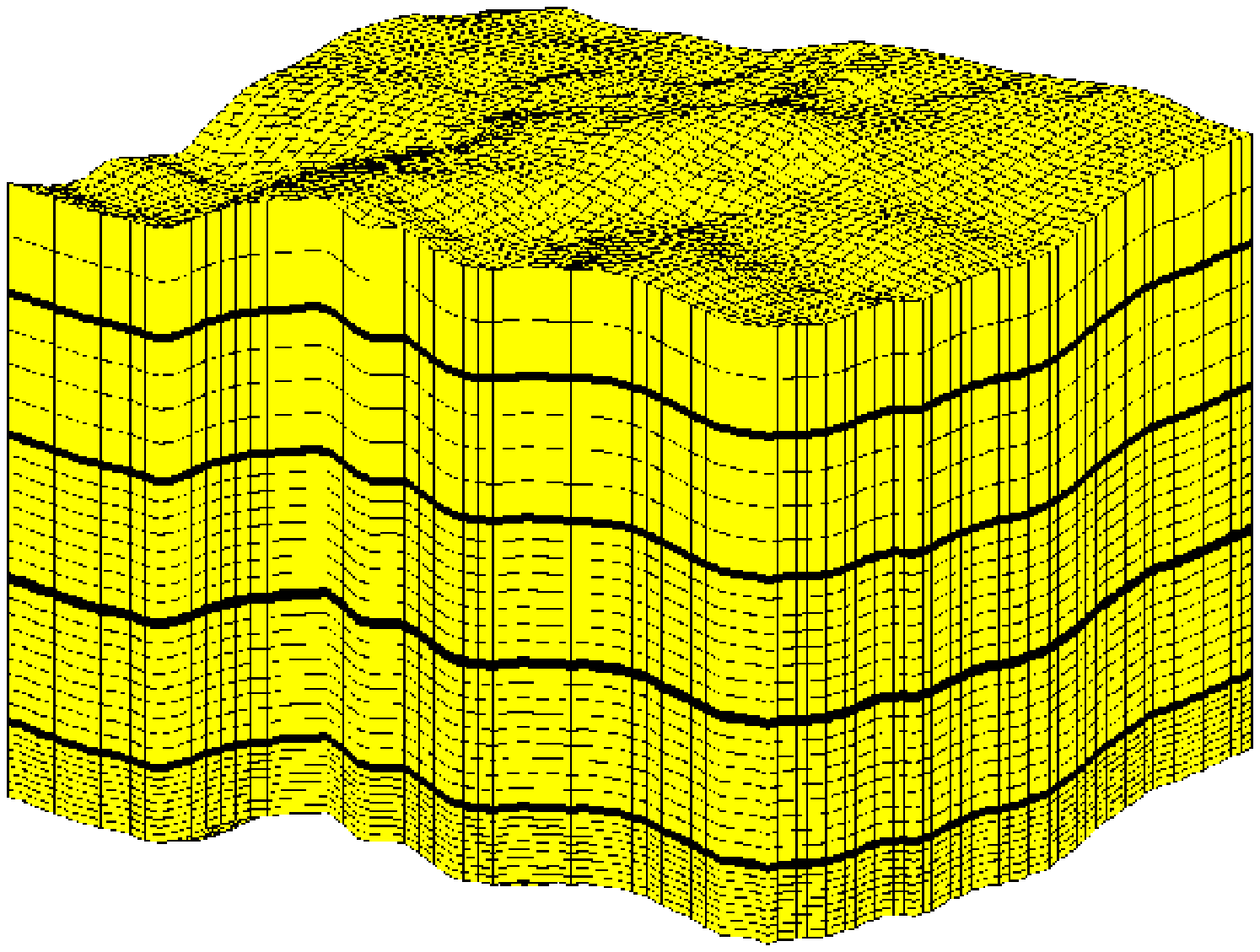  } 
\caption{A multiblock grid in a geological formation.}
\label{example1}
\end{figure}

The geological formation is shown in figure (\ref{example1}) is divided into nine layers based on the medium property. Four of these nine layers are highly permeable thus these layers are densly meshed, as shown in the figure \ref{example1}.
\section*{ACKNOWLEDGEMENTS}
We thank Ivar Aavatsmark for providing useful comments, and Many L. Buddle and David R. Wood for correcting the manuscript.

\end{spacing}
\end{document}